.
\font\piccolissimo=cmr5.
\font\script=eusm10.
\font\sets=msbm10.
\font\stampatello=cmcsc10.
\font\symbols=msam10.

\def\1{{\bf 1}}
\def\sgn{{\rm sgn}}

\def\msum{{\sum_{m\sim {N\over {\ell}}}}^{\ast}}
\def\square{\hbox{\vrule\vbox{\hrule\phantom{s}\hrule}\vrule}}
\def\defineq{\buildrel{def}\over{=}}

\def\doublesum{\mathop{\sum\sum}}

\def\N{\hbox{\sets N}}
\def\R{\hbox{\sets R}}
\def\Z{\hbox{\sets Z}}
\def\Corr{\hbox{\script C}}
\def\EssBdd{\hbox{\symbols n}\,}

\par
\noindent
\centerline{\bf ON THE SYMMETRY OF ARITHMETICAL FUNCTIONS}
\centerline{\bf IN ALMOST ALL SHORT INTERVALS, IV}
\medskip
\centerline{by}
\smallskip
\centerline{G.Coppola}
\bigskip
{
\font\eightrm=cmr8
\eightrm {
\par
{\bf Abstract.} We study the arithmetic (real) function $f=g\ast \1$, with $g$ \lq \lq essentially bounded\rq \rq \thinspace and supported over the integers of $[1,Q]$. In particular, we obtain non-trivial bounds, through $f$ \lq \lq correlations\rq \rq, for the \lq \lq Selberg integral\rq \rq \thinspace and the \lq \lq symmetry integral\rq \rq \thinspace of $f$ in almost all short intervals $[x-h,x+h]$, $N\le x\le 2N$, beyond the \lq \lq classical\rq \rq \thinspace level, up to a very high level of distribution (for $h$ not too small). This time we go beyond Large Sieve inequality [C]. Precisely, our method applies Weil bound for Kloosterman sums. 
}
\footnote{}{\par \noindent {\it Mathematics Subject Classification} $(2000) : 11N37, 11N25.$}
}
\bigskip
\par
\noindent {\bf 1. Introduction and statement of the results.}
\smallskip
\par
\noindent
This paper pursues the study of the symmetry integral (and, also, of the Selberg integral) of arithmetical functions, started in our paper [C], making a substantial progress on the level of distribution of the (sieve) arithmetical functions involved.
\par
We will study {\stampatello arithmetical functions} $f:\N \rightarrow \R$ in almost all short intervals; here, \lq \lq {\stampatello almost all}\rq \rq \thinspace $[x-h,x+h]$ \thinspace with \thinspace $N<x\le 2N$ \thinspace means all, except for at most $o(N)$ of them; where $[x-h,x+h]$ \lq \lq {\stampatello short interval}\rq \rq \thinspace means $h=h(N)$ (for $N<x\le 2N$, see above), with (as $N\to \infty$) $h\to \infty$ and $h=o(N)$. Here we'll follow Landau's notation: $F(N)=o(G(N))$ stands for (when $N\to \infty$) $F(N)/G(N)\to 0$ and $F(N)={\cal O}(G(N))$ (or the equivalent Vinogradov notation $F(N)\ll G(N)$) abbreviates (again, $N\to \infty$) $|F(N)|/G(N)$ is bounded (above, by a positive constant, called the ${\cal O}-$constant or the $\ll-$constant); $F(N)\ll_{A,B,C} G(N)$ (or $F(N)={\cal O}_{A,B,C}(G(N))$) indicates the dependence of the $\ll-$ (or ${\cal O}-$)constant on these parameters. As usual, we will indicate the greatest common divisor (GCD) of $a$ and $b$ as: $(a,b)$.
\par
We will use two important quantities (describing $f$ or $h$, respectively), namely the ({\stampatello distribution}) {\stampatello level} $\lambda$ of $f$ and the {\stampatello width} $\theta$ of the short interval. Writing $f(n)=\sum_{d|n}g(d)$, we will let $\lambda \defineq (\log Q)/(\log N)$, whenever $g(q)=0$ $\forall q>Q$. And, given $h=h(N)$ the length of our short interval(s), set $\theta \defineq (\log h)/(\log N)$. In this paper the width will be assumed to be $0<\theta<1$ (but \lq \lq far\rq \rq \thinspace from edges).
\par
As usual, we'll write $\1(n)=1$, i.e. the characteristic function of natural numbers (if no subscripts; while we write $\1_{\wp}$ to mean $1$ if $\wp$ is true, $0$ otherwise), and the above becomes $f=g\ast \1$, where $\ast$ is Dirichlet product [T]; from M\"obius inversion, see [T], $g=f\ast \mu$ (where $\mu$ is M\"obius function, $\mu(1)\defineq 1$, $\mu(n)\defineq (-1)^r$ when $n$ is the product of $r$ distinct primes, $\mu(n)\defineq 0$ in all other cases, see [T]). 
\par
Some authors, like Aurel Wintner, called\thinspace $f\ast \mu$, say\thinspace $f'$, the \lq \lq Eratosthenes transform\rq \rq \thinspace of \thinspace $f$.
\par
We call $f$ {\stampatello essentially bounded} if $\forall \varepsilon > 0$ we have $f(n)\ll n^{\varepsilon}$ (as $n\to \infty$); for example, the number of (positive) divisors of $n$, say $d(n)\defineq \sum_{d|n}1$, is such (like many other interesting arithmetical functions). Sometimes (expecially during the proofs of the Lemma and of our Theorem) we'll abbreviate $A\EssBdd B$ to mean : $\forall \varepsilon >0$ $A\ll_{\varepsilon} N^{\varepsilon} B$. Also, M\"obius inversion easily allows to prove that $f\EssBdd 1$ if, and only if, $g\EssBdd 1$ (whenever $f=g\ast \1$). In the sequel $g$ support will be $[Q,2Q]$ or $[1,Q]$ (no difference, using a dissection argument: additional logarithms are absorbed into $N^{\varepsilon}$): our arguments won't be affected (for small $Q$, [C]).
\par
The quantity checking $f$ \lq \lq almost all\rq \rq -symmetry in short intervals is the {\stampatello symmetry integral}
$$
I_f(N,h)\defineq \int_{N}^{2N}\Big| \sum_{|n-x|\le h}f(n)\sgn(n-x)\Big|^2\,dx,
$$
\par
\noindent
where \enspace $\sgn(0)\defineq 0$, \thinspace $\sgn(r)\defineq r/|r|$ ($\forall r\neq 0$). \hfil (For the definition of \thinspace $J_f$, the Selberg integral of \thinspace $f$, see [C].)
\par
\noindent
In our previous paper [C] we obtained non-trivial bounds for \thinspace $I_f$ (and $J_f$), whenever the level is \enspace $\lambda < {{1+\theta}\over 2}$. We mean, for non-trivial, bounds of the kind \enspace $I_f(N,h)\ll {{Nh^2}\over {N^{\delta}}}$ (same for $J_f$), for some \thinspace $\delta>0$ (small).
\par
Here we will supersede these results, getting level \enspace $\lambda < \max\left(1-{{\theta}\over 2},{{1+\theta}\over 2}\right)$ (using previous, too).

\vfil
\eject				

\par
\noindent
As usual, $[\alpha]\defineq {\displaystyle \min_{n\in \Z,n\le \alpha} }n$ \thinspace will be the {\stampatello integer part} of $\alpha \in \R$ and $\{ \alpha\}\defineq \alpha - [\alpha]$ its {\stampatello fractional part}.
\bigskip
\par
\noindent
We give, now, our main result.
\medskip
\par
\noindent {\stampatello Theorem.} {\it Fix } $\theta_0 \in ]0,1/2[$. {\it Let} $N,h\in \N$, {\it where} $N^{\theta}\ll h\ll N^{\theta}$, {\it if } $N\to \infty$, {\it with } $\theta_0 < \theta < 1-\theta_0$. {\it Assume } $f=g\ast \1$ {\it real and essentially bounded with, say, } $g(q)=0$, $\forall q>Q$, $Q=o({N\over {\sqrt h}})$. {\it Then}
$$
I_f(N,h)
\EssBdd Nh+h^3+Nh^2{h\over Q}+Nh^2\left({{Q^2h}\over {N^2}}\right)^{1/5}.
$$
\medskip
\par
\noindent
(We point out that a result similar to our Theorem holds for the Selberg integral of $f$, too.)
\medskip
\par
\noindent
We explicitly remark that an important limitation to obtain non-trivial bounds comes from $h=o(Q)$, but can be avoided (namely, $Q\ll h$ gives $\EssBdd h^3$, since $\lambda<{{1+\theta}\over 2}$ from quoted LS-type result of [C]).
\medskip
\par
\noindent
As regards the level, it comes from our Theorem the quoted above $\lambda<1-\theta/2$.

\bigskip

\par
\noindent
The paper is organized as follows:
\item{$\star$}
in section 2 we prove the (exponential sums) Lemma for the proof of our Theorem,
\item{$\star$}
that is given in section 3.

\vfil
\eject				

\par
\noindent {\bf 2. A trigonometric Lemma.}
\smallskip
\par
\noindent {\bf Lemma.} {\it Let } $a,j,t\in \N$ \thinspace {\it with } $(a,t)=1$ {\it and } $a=o(t)$, {\it as } $a\to \infty$, $t\to \infty$ ($a,t\ll N$, {\it if } $N\to \infty$). {\it Then}
$$
\sum_{0<|k|\le {t\over 2}}\left[ \cos{{2\pi k\overline{t}}\over a} \left( \cos{{2\pi k}\over {at}}-1\right)+\sin{{2\pi k\overline{t}}\over a}\sin{{2\pi k}\over {at}}\right]\cos{{2\pi jk}\over t} = \sum_{r\in \Z_a}\cos{{2\pi r\overline{t}}\over a}\Sigma_1 + \sum_{r\in \Z_a}\sin{{2\pi r\overline{t}}\over a}\Sigma_2,
$$
\par
\noindent
{\it defining the reciprocal } $\overline{a}a\equiv 1(\bmod\,\,t)$, {\it we have } $\sum_r \cos{{2\pi r\overline{t}}\over a}\Sigma_1(0)\ll {1\over a}$ {\it and } $\sum_r \sin{{2\pi r\overline{t}}\over a}\Sigma_2(0)\ll 1$, {\it where}
$$
\Sigma_1 = \Sigma_1(j) \defineq 2\sum_{{k\le t/2}\atop {k\equiv r(a)}}\left( \cos{{2\pi k}\over {at}}-1\right)\cos{{2\pi jk}\over t} 
\enspace \Rightarrow \enspace 
{1\over t}\sum_{0\le |j|\le {t\over 2}}\left|\sum_{0\le |r|\le {a\over 2}}\cos{{2\pi r\overline{t}}\over a}\Sigma_1(j)\right| \EssBdd {1\over a},
\leqno{(1)}
$$
\par
\noindent
({\it we abbreviate } $k\equiv r(a)$ {\it for } $k\equiv r(\bmod\,\,a)$, {\it hereon}), {\it indicating } $M(j)\defineq -2\sin{{2\pi jr}\over t}{\displaystyle \sum_{J\le {t\over {2a}}} }\sin{{2\pi J}\over t}\sin{{2\pi jaJ}\over t}$, 
$$
\Sigma_2 = \Sigma_2(j) \defineq 2\sum_{{k\le t/2}\atop {k\equiv r(a)}}\sin{{2\pi k}\over {at}}\cos{{2\pi jk}\over t}
\enspace \Rightarrow \enspace 
{1\over t}\sum_{0\le |j|\le {t\over 2}}\left|\sum_{0\le |r|\le {a\over 2}}\sin{{2\pi r\overline{t}}\over a}\left(\Sigma_2(j)\!-\!M(j)\right)\right|\EssBdd {a\over t}+{1\over a}.
\leqno{(2)}
$$
\par
\noindent {\it Proof.}$\!$ Use \thinspace $k=Ja+r$ \thinspace for (both $J=0$, see the remainder, and) $J\ge 1$:
$$
\Sigma_1 = 2\sum_{J\le {t\over {2a}}-{r\over a}}\left(\cos\left({{2\pi J}\over t}+{{2\pi r}\over {at}}\right)-1\right)\cos{{2\pi j(Ja+r)}\over t}
+ {\cal O}\left( {1\over {t^2}}\right) 
$$
$$
\enspace 
= 2\sum_{J\le {t\over {2a}}}\left( \cos{{2\pi J}\over t} - 1\right)\cos{{2\pi j(Ja+r)}\over t}+{\cal O}\left( {1\over {t^2}}+{1\over {at}}+{1\over {a^2}}\right)
$$
\par
where, due to $a=o(t)$, these remainders contribute $\EssBdd {1\over a}$ to $(1)$. Simplifying $r-$odd terms, we get
$$
\sum_{r\in \Z_a}\cos{{2\pi r\overline{t}}\over a}\left(2\sum_{J\le {t\over {2a}}}\left( \cos{{2\pi J}\over t} - 1\right)\cos{{2\pi j(Ja+r)}\over t}\right) 
$$
$$
= \sum_{r\in \Z_a}\cos{{2\pi r\overline{t}}\over a}\cos{{2\pi jr}\over t}\left(2\sum_{J\le {t\over {2a}}}\left( \cos{{2\pi J}\over t} - 1\right)\cos{{2\pi jaJ}\over t}\right)\ll {1\over a}\left(\1_{j\neq 0}{1\over {\left\Vert {{ja}\over t}\right\Vert}}\right),
$$
\par
hereon $\Vert \alpha \Vert \defineq {\displaystyle \min_{n\in \Z} }|\alpha -n|$, applying partial summation [T] to ${\displaystyle \sum_J }e(J\alpha)\ll {1\over {\left\Vert \alpha \right\Vert}}$, see [D,chap.26] ($\alpha={{ja}\over t}$);
$$
{1\over t}\sum_{j\neq 0}{1\over {\left\Vert {{ja}\over t}\right\Vert}}=\sum_{n\le t/2}{1\over n}\sum_{j\equiv \pm \overline{a}n(t)}1=2\sum_{n\le {t\over 2}}{1\over n}\ll \log(t+2)\EssBdd 1.
$$
\par
Then, $(1)$ is proved. 
\par
In order to prove $(2)$, we first let \thinspace $k=Ja+r$ \thinspace (again, for $J=0$, see the remainder): 
$$
\Sigma_2(j) = {\cal O}\left({1\over t}\right) + 2\sum_{J\le {t\over {2a}}-{r\over a}}\sin\left({{2\pi J}\over t}+{{2\pi r}\over {at}}\right)\cos\left({{2\pi jaJ}\over t}+{{2\pi jr}\over t}\right)
$$
\par
We exclude \thinspace $j\equiv \pm \overline{a}(t)$ \thinspace in the following.
\par
We sweep out, this time, after a small correction like in the previous Lemma, the $r-$even terms, 
$$
\sum_{0<|r|\le a/2}\sin{{2\pi r\overline{t}}\over a}(\Sigma_2 -M(j))
= 2\sum_{0<|r|\le a/2}\sin{{2\pi r\overline{t}}\over a}\sum_{J\le {t\over {2a}}}\sin\left({{2\pi J}\over t}+{{2\pi r}\over {at}}\right)\cos\left({{2\pi jaJ}\over t}+{{2\pi jr}\over t}\right) 
$$

\vfil
\eject				

$$
+ {\cal O}\left({a\over t}\right) - 2\sum_{r\le a/2}\sin{{2\pi r\overline{t}}\over a}\sum_{{t\over {2a}}<J\le {t\over {2a}}+{r\over a}}\sin\left({{2\pi J}\over t}-{{2\pi r}\over {at}}\right)\cos\left({{2\pi jaJ}\over t}-{{2\pi jr}\over t}\right) 
$$
$$
- 2\sum_{r\le a/2}\sin{{2\pi r\overline{t}}\over a}\sum_{{t\over {2a}}-{r\over a}<J\le {t\over {2a}}}\sin\left({{2\pi J}\over t}+{{2\pi r}\over {at}}\right)\cos\left({{2\pi jaJ}\over t}+{{2\pi jr}\over t}\right) 
$$
$$
+ 2\sum_{0<|r|\le a/2}\sin{{2\pi r\overline{t}}\over a}\sin{{2\pi jr}\over t}\sum_{J\le {t\over {2a}}}\sin{{2\pi J}\over t}\sin{{2\pi jaJ}\over t} 
$$
$$
= {\cal O}\left({a\over t}\right) 
- 2\left(\sum_{0<|r|\le a/2}\sin{{2\pi r\overline{t}}\over a}\left(\cos{{2\pi r}\over {at}}-1\right)\sin{{2\pi jr}\over t}\right)\left(\sum_{J\le {t\over {2a}}}\sin{{2\pi J}\over t}\sin{{2\pi jaJ}\over t}\right) 
$$
$$
+ 2\left(\sum_{0<|r|\le a/2}\sin{{2\pi r\overline{t}}\over a}\sin{{2\pi r}\over {at}}\cos{{2\pi jr}\over t}\right)\left(\sum_{J\le {t\over {2a}}}\cos{{2\pi J}\over t}\cos{{2\pi jaJ}\over t}\right) 
$$
$$
- 2\sum_{r\le a/2}\sin{{2\pi r\overline{t}}\over a}\sum_{{t\over {2a}}<J\le {t\over {2a}}+{r\over a}}\sin{{2\pi J}\over t}\cos\left({{2\pi jaJ}\over t}-{{2\pi jr}\over t}\right) + {\cal O}\left({a\over t}\right) 
$$
$$
- 2\sum_{r\le a/2}\sin{{2\pi r\overline{t}}\over a}\sum_{{t\over {2a}}-{r\over a}<J\le {t\over {2a}}}\sin{{2\pi J}\over t}\cos\left({{2\pi jaJ}\over t}+{{2\pi jr}\over t}\right) + {\cal O}\left({a\over t}\right); 
$$
\par
now on, we'll ignore all remainders giving a contribution already into $(2)$: the first product here is $\ll {1\over t}$.
\par
Then, from the quoted [D] above, the second (\lq \lq non-sporadic\rq \rq) term is ($\forall j\not \equiv \pm \overline{a}(t)$,here)
$$
\thinspace \enspace \thinspace \thinspace 
{\cal O}\left( {a\over t}\right) \sum_{J\le {t\over {2a}}}\cos{{2\pi J}\over t}\cos{{2\pi jaJ}\over t} 
= {\cal O}\left( {a\over t}\right) \sum_{J\le {t\over {2a}}}\cos{{2\pi jaJ}\over t} + {\cal O}\left( {1\over {a^2}}\right) 
\ll \1_{j=0}+{a\over t}\left(\1_{j\neq 0}{1\over {\left \Vert {{ja}\over t}\right \Vert}}\right)+{1\over {a^2}}
$$
\par
Since we have \thinspace $\sin{{2\pi J}\over t}={{\pi}\over a}+{\cal O}({1\over t})$ \thinspace into sporadic terms, these are:
$$
\quad \enspace \thinspace 
- {{2\pi}\over a}\sum_{r\le a/2}\sin{{2\pi r\overline{t}}\over a}\left(\sum_{{t\over {2a}}<J\le {t\over {2a}}+{r\over a}}\cos\left({{2\pi jaJ}\over t}-{{2\pi jr}\over t}\right) + \sum_{{t\over {2a}}-{r\over a}<J\le {t\over {2a}}}\cos\left({{2\pi jaJ}\over t}+{{2\pi jr}\over t}\right)\right) 
$$
\par
($\Rightarrow $ \thinspace $\sum_r \sin{{2\pi r\overline{t}}\over a}\Sigma_2(0)\ll 1$), i.e., applying {\stampatello reciprocity} ${{\overline{t}}\over a}\equiv -{{\overline{a}}\over t} + {1\over {at}}(\bmod\,\,1)$ and expanding cosines, 
$$
- {{2\pi}\over a}\cos{{2\pi ja}\over t}\!\left[ {t\over {2a}}+{1\over 2}\right]\sum_{r\le a/2}\left( \left[ {t\over {2a}}+{r\over a}\right] - \left[ {t\over {2a}}\right]\right) \sin{{2\pi r\overline{a}}\over t}\cos{{2\pi jr}\over t} 
$$
$$
- {{2\pi}\over a}\sin{{2\pi ja}\over t}\!\left[ {t\over {2a}}+{1\over 2}\right]\sum_{r\le a/2}\left( \left[ {t\over {2a}}+{r\over a}\right] - \left[ {t\over {2a}}\right]\right)\sin{{2\pi r\overline{a}}\over t}\sin{{2\pi jr}\over t}
$$
$$
- {{2\pi}\over a}\cos{{2\pi ja}\over t}\!\left[ {t\over {2a}}\right]\sum_{r\le a/2}\left( \left[ {t\over {2a}}\right] - \left[ {t\over {2a}}-{r\over a}\right]\right) \sin{{2\pi r\overline{a}}\over t}\cos{{2\pi jr}\over t} 
$$
$$
\thinspace \quad \enspace \thinspace 
+ {{2\pi}\over a} \sin{{2\pi ja}\over t}\!\left[ {t\over {2a}}\right]\sum_{r\le a/2}\left( \left[ {t\over {2a}}\right] - \left[ {t\over {2a}}-{r\over a}\right]\right)\sin{{2\pi r\overline{a}}\over t}\sin{{2\pi jr}\over t} + {\cal O}\left( {1\over t}\right) 
$$
$$
\thinspace \quad \thinspace 
\ll {1\over a}\!\left(\left|\sum_{R_1<r\le R_2}\!e_t((\overline{a}\!+\!j)r)\right|\!+\!\left|\sum_{R_1<r\le R_2}\!e_t((\overline{a}\!-\!j)r)\right|\!+\!\left|\sum_{R_3<r\le R_4}\!e_t((\overline{a}\!+\!j)r)\right|\!+\!\left|\sum_{R_3<r\le R_4}\!e_t((\overline{a}\!-\!j)r)\right|\right) 
\!\!+\!{1\over t} 
$$
$$
\ll {a\over t} + {1\over a}\left({1\over {\left \Vert {{j+\overline{a}}\over t}\right \Vert}}+{1\over {\left \Vert {{j-\overline{a}}\over t}\right \Vert}}\right).\enspace \square
$$

\vfil
\eject

\par				
\noindent {\bf 3. Proof of the Theorem.}
\smallskip
\par
\noindent
A kind of \lq \lq elementary dispersion\rq \rq ([C, Lemma 1]) gives (writing $n\sim N$ in sums for $N<n\le 2N$), defining
$$
W(a)\defineq \1_{[-2h,2h]}(a)\max(2h-3|a|,|a|-2h) \enspace \quad \hbox{\stampatello and} \quad \enspace \Corr_f(a)\defineq \sum_{n\sim N}f(n)f(n-a),
$$
$$
I_f(N,h)=\sum_{a\neq 0}W(a)\Corr_f(a)+{\cal O}\left( N^{\varepsilon}\left( Nh+h^3\right)\right).
$$
\par
\noindent
(We'll refer to the first kind of remainders as the \lq \lq {\stampatello diagonal-type}\rq \rq, while the others will be the \lq \lq {\stampatello tails}\rq \rq.)
\par
\noindent
Then, after using $f(n-a)=\sum_{q|n-a}g(q)$, we organize the $q-$sum w.r.t. $(a,q)={\rm GCD}(a,q)$, say:
$$
\sum_{a\neq 0}W(a)\Corr_f(a) = \sum_{a\neq 0}W(a)\sum_{\ell | a}\sum_{(q,a)=\ell}g(q)\sum_{{n\sim N}\atop {n\equiv a(q)}}f(n)
$$
\par
\noindent
which, changing variables and using $|a|\le 2h$ (otherwise $W=0$), is 
$$
\sum_{a\neq 0}W(a)\Corr_f(a) = \sum_{\ell \le 2h}\sum_{a\neq 0}W(a\ell)\sum_{(q,a)=1}g(\ell q)\sum_{{m\sim {N\over {\ell}}}\atop {m\equiv a(q)}}f(\ell m), 
$$
\par
\noindent
then (using $Q\ll N$ from $g$ definition and \thinspace $W(a)\ll h$ \thinspace {\stampatello uniformly} on $a$) we \lq \lq cut\rq \rq \thinspace the divisors $\ell$ at $G=o(h)$:
$$
\sum_{a\neq 0}W(a)\Corr_f(a) = \sum_{\ell \le G}\sum_{a\neq 0}W(a\ell)\sum_{(q,a)=1}g(\ell q)\sum_{{m\sim {N\over {\ell}}}\atop {m\equiv a(q)}}f(\ell m) + {\cal O}\left( N^{\varepsilon} {{Nh^2}\over G}\right);
$$
\par
\noindent
also, we may assume $|a|>h/G$, otherwise the error is again the same (using $h=o(N)$, here):
$$
\sum_{a\neq 0}W(a)\Corr_f(a) = \sum_{\ell \le G}\sum_{|a|>{h\over G}}W(a\ell)\sum_{(q,a)=1}g(\ell q){\sum_{{m\sim {N\over {\ell}}}\atop {m\equiv a(q)}}}^{\ast}f(\ell m) + {\cal O}\left( N^{\varepsilon} \left({{Nh^2}\over G}\right)\right).
$$
\par
\noindent
Hereon we'll abbreviate the condition \lq \lq $(m,q)=1$\rq \rq \thinspace with a \lq \lq $\ast$\rq \rq \thinspace in the $m-$sum.
\par
\noindent
(We use coprimality of $a$ and $q$ to insert a limitation which is immaterial: $m\equiv a(q)$ $\Rightarrow$ $(m,q)=1$.)\par
\noindent
Since remainders are negligible, we write \thinspace $\sim$ \thinspace to leave them. \hfil (No confusion arises with \thinspace $\sim$ \thinspace in summations !) 
\par
\noindent
Thus
$$
I_f(N,h)\sim \sum_{\ell \le G}\sum_{|a|>{h\over G}}W(a\ell)\sum_{(q,a)=1}g(\ell q){\sum_{{m\sim {N\over {\ell}}}\atop {m\equiv a(q)}}}^{\ast}f(\ell m) 
$$
\par
\noindent
and, applying the orthogonality [V] of additive characters \thinspace $e_q(ja)\defineq e({{ja}\over q})$, where $e(\alpha)\defineq e^{2\pi i\alpha}$, this is
$$
I_f(N,h)\sim \sum_{\ell \le G}\sum_{|a|>{h\over G}}W(a\ell)\sum_{(q,a)=1}{{g(\ell q)}\over q}\sum_{0<|j|\le q/2}e_q(ja)\msum f(\ell m)e_q(-jm) 
$$
\par
\noindent
since $j=0$ terms give diagonal-type remainders, due to: ${\displaystyle \sum_{a\neq 0} }W(a\ell)\ll h$ [C,Lemma 4,(1)]. \thinspace Being $W$ {\stampatello even}
$$
I_f(N,h)
\sim \sum_{\ell \le G}\sum_{|a|>{h\over G}}W(a\ell)\sum_{(q,a)=1}{{g(\ell q)}\over q}\sum_{0<|j|\le {q\over 2}}\cos{{2\pi ja}\over q}\msum f(\ell m)\cos{{2\pi jm}\over q} 
$$
$$
\sim 2\sum_{\ell \le G}\sum_{{h\over G}<a\le {{2h}\over {\ell}}}W(a\ell)\sum_{(q,a)=1}{{g(\ell q)}\over q}\sum_{0<|j|\le {q\over 2}}\cos{{2\pi ja}\over q}\msum f(\ell m)\cos{{2\pi jm}\over q} 
$$
\par
\noindent			
and, now, we exploit the coprimality of $a$ with $q$, changing variable:
$$
I_f(N,h)\sim 2\sum_{\ell \le G}\sum_{{h\over G}<a\le {{2h}\over {\ell}}}W(a\ell)\sum_{(q,a)=1}{{g(\ell q)}\over q}\sum_{0<|j|\le q/2}\cos{{2\pi j}\over q}\msum f(\ell m)\cos{{2\pi jm\overline{a}}\over q} 
$$
\par
\noindent
We rearrange, then, the terms, depending on residue classes modulo $q$:
$$
I_f(N,h)\sim 2\sum_{\ell \le G}\sum_{{h\over G}<a\le {{2h}\over {\ell}}}W(a\ell)\sum_{(q,a)=1}{{g(\ell q)}\over q}\doublesum_{{0<|j|,|k|\le q/2}\atop {(j,q)=(k,q)}}\cos{{2\pi j}\over q}\cos{{2\pi k\overline{a}}\over q}{\sum_{{m\sim {N\over {\ell}}}\atop {jm\equiv k(q)}}}^{\ast}f(\ell m)
$$
\par
\noindent
Therefore, applying {\stampatello reciprocity} \thinspace ${{\overline{a}}\over q} \equiv -{{\overline{q}}\over a}+{1\over {aq}}(\bmod \,\, 1)$, \thinspace we separate into two sums, say $I_f \sim S_1 + S_2$ :
$$
S_1 \defineq 2\sum_{\ell \le G}\sum_{{h\over G}<a\le {{2h}\over {\ell}}}W(a\ell)\sum_{(q,a)=1}{{g(\ell q)}\over q}\doublesum_{0<|j|,|k|\le q/2}\cos{{2\pi j}\over q}\cos{{2\pi k\overline{q}}\over a}{\sum_{{m\sim {N\over {\ell}}}\atop {jm\equiv k(q)}}}^{\ast}f(\ell m)
$$
\par
\noindent
and
$$
S_2 \defineq 2\sum_{\ell \le G}\sum_{{h\over G}<a\le {{2h}\over {\ell}}}W(a\ell)\sum_{(q,a)=1}{{g(\ell q)}\over q}\Sigma,
$$
\par
\noindent
where, say, 
$$
\Sigma \defineq \doublesum_{{0<|j|,|k|\le q/2}\atop {(j,q)=(k,q)}}\cos{{2\pi j}\over q}
\left( \cos{{2\pi k\overline{q}}\over a} \left( \cos{{2\pi k}\over {aq}}-1\right)+\sin{{2\pi k\overline{q}}\over a}\sin{{2\pi k}\over {aq}}\right)
{\sum_{{m\sim {N\over {\ell}}}\atop {jm\equiv k(q)}}}^{\ast}f(\ell m)
$$
\par
\noindent
We'll treat first $S_2$, as it contains \lq \lq remainder terms\rq \rq, here. Collect all $j$ with same $(j,q)$ and change $j$, $k$:
$$
\Sigma = \sum_{{t|q}\atop {t<q}}\doublesum_{{0<|j|,|k|\le {q\over {2t}}\atop {(j,{q\over t})=(k,{q\over t})=1}}}\cos{{2\pi j}\over {q/t}}\left( \cos{{2\pi k\overline{q/t}}\over a} \left( \cos{{2\pi k}\over {a(q/t)}}-1\right)+\sin{{2\pi k\overline{q/t}}\over a}\sin{{2\pi k}\over {a(q/t)}}\right){\sum_{{m\sim {N\over {\ell}}}\atop {m\equiv \overline{j}k(q/t)}}}^{\ast}f(\ell m)
\buildrel{{\rm F}}\over{=}
$$
$$
\buildrel{{\rm F}}\over{=}
\sum_{{t|q}\atop {t>1}}\sum_{{0<|k|\le {t\over 2}}\atop {(k,t)=1}}\left( \cos{{2\pi k\overline{t}}\over a} \left( \cos{{2\pi k}\over {at}}-1\right)+\sin{{2\pi k\overline{t}}\over a}\sin{{2\pi k}\over {at}}\right)\sum_{{0<|j|\le {t\over 2}}\atop {(j,t)=1}}\cos{{2\pi j}\over t}{\sum_{{m\sim {N\over {\ell}}}\atop {m\equiv \overline{j}k(t)}}}^{\ast}f(\ell m)
$$
\par
\noindent
where \lq \lq {\stampatello f}\rq \rq \thinspace means we \lq \lq flipped\rq \rq \thinspace the divisors $t$. Our inner sum over $j,m$ is, again through orthogonality, 
$$
\sum_{{0<|j|\le {t\over 2}}\atop {(j,t)=1}}\cos{{2\pi j}\over t}{\sum_{{m\sim {N\over {\ell}}}\atop {m\equiv \overline{j}k(t)}}}^{\ast}f(\ell m) 
= {1\over t}\sum_{{0<|j|\le {t\over 2}}\atop {(j,t)=1}}\cos{{2\pi j}\over t}\sum_{0\le |s|\le t/2}e_t(sk\overline{j})\msum f(\ell m)e_t(-sm)
$$
$$
= {{\mu(t)}\over t}\msum f(\ell m) + {1\over {2t}}\sum_{0<|s|\le t/2}\left[ S(1,sk;t)+S(-1,sk;t)\right]\msum f(\ell m)e_t(-sm),
$$
\par
\noindent
where $\mu$ is M\"obius function (see above, in the introduction),
$$
c_t(n)\defineq \sum_{j\le t,(j,t)=1}e_t(jn)=\varphi(t){{\mu(t/(t,n))}\over {\varphi(t/(t,n))}}=\sum_{{0<|j|\le {t\over 2}}\atop {(j,t)=1}}\cos{{2\pi jn}\over t} \quad \hbox{\rm is}\enspace \hbox{\stampatello Ramanujan\thinspace sum}
$$
\par				
\noindent
($n=1$ gives $\mu(t)$ here), see quoted [D], with Euler function $\varphi(n)\defineq \sum_{j\le n,(j,n)=1}1$ and
$$
S(a,b;c)\defineq \sum_{{j\le c,(j,c)=1}\atop {j\overline{j}\equiv 1(c)}}e_c(ja+\overline{j}b) \quad \hbox{\rm is}\enspace \hbox{\stampatello Kloosterman\thinspace sum}
$$
\par
\noindent
that can be estimated using Weil-Estermann bound, see [I-K, chap.11],
$$
S(a,b;c)\ll d(t)\sqrt{(a,b,c)}\sqrt{c} \quad ( \Rightarrow \enspace \left[S(1,ks;t)+S(-1,ks;t)\right]\EssBdd \sqrt{t}, \hbox{\rm here}).
$$
\par
\noindent
We'll use this bound in a while. First, see that (recall $s\not \equiv 0(t)$ and use $\1_{(m,q)=1}=\sum_{d|q,d|m}\mu(d)$, [T])
$$
\msum f(\ell m)e_t(-sm) = \sum_{q'}g(q')\sum_{{n\sim N}\atop {{n\equiv 0(q'),n\equiv 0(\ell)}\atop {({n\over {\ell}},q)=1}}}e_t(-s(n/\ell)) 
$$
$$
= \sum_{\ell' | \ell} \sum_{{q''}\atop {(q'',{{\ell}\over {\ell'}}q)=1}}g(\ell' q'')\sum_{{m\sim N/(\ell q'')}\atop {(m,q)=1}}e_t(-sq''m) 
= \sum_{\ell' | \ell} \sum_{{q''}\atop {(q'',{{\ell}\over {\ell'}}q)=1}}g(\ell' q'')\sum_{d|q}\mu(d)\sum_{m\sim N/(\ell q''d)}e_t(-sq''dm)
$$
$$
= \sum_{\ell' | \ell} \sum_{{q''}\atop {(q'',{{\ell}\over {\ell'}}q)=1}}g(\ell' q'')\sum_{{d|q}\atop {sq''d\equiv 0(t)}}\mu(d)\left({N\over {\ell q''d}}+{\cal O}(1)\right) + {\cal O}\left(N^{\varepsilon}\sum_{\ell' | \ell} \sum_{{q''\sim {Q\over {\ell'}}}\atop {(q'',{{\ell}\over {\ell'}}q)=1}}\sum_{{d|q}\atop {sq''d\not \equiv 0(t)}}{1\over \left \Vert {sq''d\over t}\right \Vert}\right)
$$
\par
\noindent
(we used, when $\Vert \beta \Vert \neq 0$, the classical bound \thinspace $\sum_{N_1\le n\le N_2}e(n\beta)\ll 1/\Vert \beta \Vert$, \thinspace see quoted [D]) gives
$$
\sum_{(q,a)=1}{{g(\ell q)}\over q}\Sigma = \sum_{{1<t\le {{2Q}\over {\ell}}}\atop {(t,a)=1}}\left( \sum_{(q,a)=1}{{g(\ell tq)}\over q}\right){1\over t}\sum_{{0<|k|\le {t\over 2}}\atop {(k,t)=1}}\left( \cos{{2\pi k\overline{t}}\over a} \left( \cos{{2\pi k}\over {at}}-1\right)+\sin{{2\pi k\overline{t}}\over a}\sin{{2\pi k}\over {at}}\right)\times
\leqno{(*)}
$$
$$
\times \left( {{\mu(t)}\over t}\sum_{{m\sim {N\over {\ell}}}\atop {(m,qt)=1}}f(\ell m) + {1\over {2t}}\sum_{0<|s|\le t/2}\left[ S(1,sk;t)+S(-1,sk;t)\right]\left[ \sum_{\ell' | \ell} \sum_{{q''}\atop {(q'',{{\ell}\over {\ell'}}qt)=1}}g(\ell' q'')\sum_{{d|qt}\atop {sq''d\equiv 0(t)}}\mu(d)\times
\right. \right.
$$
$$
\left. \left. 
\times \left({N\over {\ell q''d}} + {\cal O}(1)\right) + {\cal O}\left(N^{\varepsilon}\sum_{\ell' | \ell} \sum_{{q''\sim {Q\over {\ell'}}}\atop {(q'',{{\ell}\over {\ell'}}qt)=1}}\sum_{{d|qt}\atop {sq''d\not \equiv 0(t)}}{1\over \left \Vert {sq''d\over t}\right \Vert}\right)\right]\right)
$$
\par
\noindent
which we'll distinguish into $t\le {T\over {\ell}}$ and $t>{T\over {\ell}}$; for the former
$$
\sum_{{1<t\le {T\over {\ell}}}\atop {(t,a)=1}}\left( \sum_{(q,a)=1}{{g(\ell tq)}\over q}\right){1\over t}\sum_{{0<|k|\le {t\over 2}}\atop {(k,t)=1}}\left( \cos{{2\pi k\overline{t}}\over a} \left( \cos{{2\pi k}\over {at}}-1\right)+\sin{{2\pi k\overline{t}}\over a}\sin{{2\pi k}\over {at}}\right)\times
$$
$$
\times \left( {{\mu(t)}\over t}\sum_{{m\sim {N\over {\ell}}}\atop {(m,qt)=1}}f(\ell m) + {1\over {2t}}\sum_{0<|s|\le t/2}
\left[ S(1,sk;t)+S(-1,sk;t)\right]
\left[ \sum_{\ell' | \ell} \sum_{{q''}\atop {(q'',{{\ell}\over {\ell'}}qt)=1}}g(\ell' q'')\sum_{{d|qt}\atop {sq''d\equiv 0(t)}}\mu(d){\cal O}(1) 
\right.
\right.
$$
$$				
\left.
\left.
+ {\cal O}\left(N^{\varepsilon}\sum_{\ell' | \ell} \sum_{{q''\sim {Q\over {\ell'}}}\atop {(q'',{{\ell}\over {\ell'}}qt)=1}}\sum_{{d|qt}\atop {sq''d\not \equiv 0(t)}}{1\over \left \Vert {sq''d\over t}\right \Vert}\right)\right]\right)
\EssBdd {N\over {a\ell}} + {{QT^{3/2}}\over {a\ell^{3/2}}} + {{T^{5/2}}\over {a\ell^{5/2}}},
$$
\par
\noindent
having used the trivial bound \enspace ${1\over t}{\displaystyle \sum_{{0<|k|\le {t\over 2}}\atop {(k,t)=1}} }\left( \cos{{2\pi k\overline{t}}\over a} \left( \cos{{2\pi k}\over {at}}-1\right)+\sin{{2\pi k\overline{t}}\over a}\sin{{2\pi k}\over {at}}\right)\ll {1\over a}$ \enspace and, esp., 
$$
\sum_{{1<t\le {T\over {\ell}}}\atop {(t,a)=1}}\left( \sum_{(q,a)=1}{{g(\ell tq)}\over q}\right){1\over t}\sum_{{0<|k|\le {t\over 2}}\atop {(k,t)=1}}\left( \cos{{2\pi k\overline{t}}\over a} \left( \cos{{2\pi k}\over {at}}-1\right)+\sin{{2\pi k\overline{t}}\over a}\sin{{2\pi k}\over {at}}\right)\times
$$
$$
\times {1\over {2t}}\sum_{0<|s|\le t/2}\left[ S(1,sk;t)+S(-1,sk;t)\right]
{\cal O}\left(N^{\varepsilon}\sum_{\ell' | \ell} \sum_{{q''\sim {Q\over {\ell'}}}\atop {(q'',{{\ell}\over {\ell'}}qt)=1}}\sum_{{d|qt}\atop {sq''d\not \equiv 0(t)}}{1\over \left \Vert {sq''d\over t}\right \Vert}\right)
$$
$$
\EssBdd \sum_{t\le {T\over {\ell}}}{1\over a}{1\over t}\sum_{s\le {t\over 2}}{\sqrt{t}}\sum_{\ell' | \ell}\sum_{q''\sim {Q\over {\ell'}}}\sum_{d\le {{2Q}\over {\ell}}}{1\over d}\sum_{{r\le t/2}\atop {sq''d\equiv \pm r(t)}}{t\over {r}}
\EssBdd \sum_{t\le {T\over {\ell}}}{{\sqrt{t}}\over a}\sum_{t'|t}\left({Q\over {t/t'}}+1\right)t\sum_{r\le {t\over {2t'}}}{1\over {t'r}}
\EssBdd {{QT^{3\over 2}}\over {a\ell^{3\over 2}}} + {{T^{5\over 2}}\over {a\ell^{5\over 2}}}.
$$
\par
\noindent
Here we use the Weil bound for Kloosterman sums. Last term of \thinspace $t\le T/\ell$ \thinspace is
$$
\sum_{{1<t\le {T\over {\ell}}}\atop {(t,a)=1}}\left( \sum_{(q,a)=1}{{g(\ell tq)}\over q}\right){1\over t}\sum_{{0<|k|\le {t\over 2}}\atop {(k,t)=1}}\left( \cos{{2\pi k\overline{t}}\over a} \left( \cos{{2\pi k}\over {at}}-1\right)+\sin{{2\pi k\overline{t}}\over a}\sin{{2\pi k}\over {at}}\right)
\times 
$$
$$
\times 
{1\over {2t}}\sum_{0<|s|\le t/2}\left[ S(1,sk;t)+S(-1,sk;t)\right]{N\over {\ell}}\sum_{\ell' | \ell} \sum_{{q''}\atop {(q'',{{\ell}\over {\ell'}}qt)=1}}{{g(\ell' q'')}\over {q''}}\sum_{{d|qt}\atop {sq''d\equiv 0(t)}}{{\mu(d)}\over d}\EssBdd {N\over {a\ell}}
$$
\par
\noindent
since terms after $\times$ can be treated opening the Kloosterman sum, see the above, becoming 
$$
{N\over {\ell t}}\sum_{(j,t)=1}\cos{{2\pi j}\over t}\sum_{t'|t}\sum_{\ell' | \ell} \sum_{{q''}\atop {(q'',{{\ell}\over {\ell'}}qt)=1}}{{g(\ell' q'')}\over {q''}}\sum_{{d|qt}\atop {(d,t)=t'}}{{\mu(d)}\over d}\sum_{{0<|s|\le t/2}\atop {s\equiv 0(t/t')}}e_t(k\overline{j}s),
$$
\par
\noindent
but this time we gain $t$, from $\sum_{(j,t)=1}\cos{{2\pi j}\over t}=\mu(t)$ (due to independence of the following sum from $j$):
$$
\sum_{{0<|s|\le t/2}\atop {s\equiv 0(t/t')}}e_t(k\overline{j}s)=\sum_{0<|s|\le t'/2}e_{t'}(k\overline{j}s)=\1_{k\equiv 0(t')}-{1\over {t'}}.
$$
\par
\noindent
Finally, the $t\le T/\ell$ terms of $(\ast)$ are in total $\EssBdd {N\over {a\ell}} + {{QT^{3/2}}\over {a\ell^{3/2}}} + {{T^{5/2}}\over {a\ell^{5/2}}}$; contribute to $S_2$ as $\EssBdd (Q+T)T^{3/2}h$.
\smallskip
\par
We pass to the terms $t>T/\ell$.
\smallskip
\par
\noindent
For these (hence, $t$ is \lq \lq {\stampatello large}\rq \rq), we use Lemma 1:
$$
\sum_{{t|q}\atop {t\,\hbox{\piccolissimo LARGE}}}\sum_{{0<|j|\le {t\over 2}}\atop {(j,t)=1}}\sum_{0<|k|\le {t\over 2}}\left[ \cos{{2\pi k\overline{t}}\over a} \left( \cos{{2\pi k}\over {at}}-1\right)+\sin{{2\pi k\overline{t}}\over a}\sin{{2\pi k}\over {at}}\right]\cos{{2\pi jk}\over t}{\sum_{{m\sim {N\over {\ell}}}\atop {m\equiv \overline{j}(t)}}}^{\ast} f(\ell m) 
$$
$$				
= -2\sum_{{t|q}\atop {t\,\hbox{\piccolissimo LARGE}}}\sum_{{|j|\le {t\over 2}}\atop {(j,t)=1}}\sum_{|r|\le {a\over 2}}\sin{{2\pi r\overline{t}}\over a}\sin{{2\pi jr}\over t}\sum_{J\le {t\over {2a}}}\sin{{2\pi J}\over t}\sin{{2\pi jaJ}\over t}{\sum_{{m\sim {N\over {\ell}}}\atop {jm\equiv 1(t)}}}^{\ast} f(\ell m) 
+ {\cal O}\!\left(\!N^{\varepsilon}\!\!\sum_{{t|q}\atop {t\,\hbox{\piccolissimo LARGE}}}\!\!\left(\!{1\over a}\!+\!{a\over t}\right)\!{N\over {\ell}}\!\right) 
$$
\par
\noindent
so to apply quoted [D] ${\displaystyle \sum_{N_1<n\le N_2} }e(\alpha n)\ll \min(N_2-N_1,{1\over {\left \Vert \alpha\right \Vert}})$ and partial summation (with {\stampatello reciprocity})
$$
\sum_{{0<|j|\le {t\over 2}}\atop {j\not \equiv \pm \overline{a}(t)}}\left|\sum_{|r|\le {a\over 2}}\sin{{2\pi r\overline{t}}\over a}\sin{{2\pi jr}\over t}\right|\left|\sum_{J\le {t\over {2a}}}\sin{{2\pi J}\over t}\sin{{2\pi jaJ}\over t}\right| 
\ll {1\over a}\sum_{{0<|j|\le {t\over 2}}\atop {j\not \equiv \pm \overline{a}(t)}}\min\left(a,{1\over {\left \Vert {{j\pm \overline{a}}\over t}\right \Vert}}\right) \min\left( {t\over a},{1\over {\left \Vert {{ja}\over t}\right \Vert}}\right) 
$$
$$
\ll {1\over a}\sum_{1<r\le {t\over 2}}\left( \sum_{j\equiv \pm \overline{a}r(t)}1\right)\min\left(a,{1\over {\left \Vert {{(r\pm 1)\overline{a}}\over t}\right \Vert}}\right) \min\left( {t\over a},{t\over r}\right) 
\ll {1\over a}\sum_{r<a}{t\over a}\min\left(a,{1\over {\left \Vert {{r\overline{a}}\over t}\right \Vert}}\right) + {t\over a} + {1\over a}\sum_{j\le {t\over {2a}}}{t\over j}
$$
$$
+ {1\over a}\sum_{{a<r\le {t\over 2}}\atop {r\not \equiv 0 (a)}}{t\over r}\min\left(a,{1\over {\left \Vert {{r\overline{a}}\over t}\right \Vert}}\right) 
\ll {1\over a}\sum_{0<|r|\le {a\over 2}}{t\over a}\min\left(a,{1\over {\left \Vert {{r\overline{t}}\over a}\right \Vert - {1\over {2a}}}}\right) + {1\over a}\sum_{{a<r\le {t\over 2}}\atop {r\not \equiv 0 (a)}}{t\over r}\min\left(a,{1\over {\left \Vert {{r\overline{t}}\over a}\right \Vert - {|r|\over {at}}}}\right) + {t\over a}\log t
$$
$$
\ll {1\over a}\left({t\over a}\sum_{0<|r|\le {a\over 2}}\min\left(a,{1\over {{{|r|}\over a}-{1\over {2a}}}}\right)+t\sum_{j\le {t\over {2a}}}{1\over {ja}}\sum_{0<|r-ja|\le {a\over 2}}\min\left(a,{1\over {{{|r-ja|}\over a}-{1\over {2a}}}}\right)+t\log t\right) 
\EssBdd {t\over a}.
$$
\par
\noindent
Here we may assume $j\not \equiv \pm \overline{a}(t)$, since the corresponding terms give a contribute already in those before: 
$$
\EssBdd \sum_{t|q,t\,\hbox{\piccolissimo LARGE}}\left|\sum_{|r|\le {a\over 2}}\sin{{2\pi r\overline{t}}\over a}\sin{{2\pi \overline{a}r}\over t}\right|\left|\sum_{J\le {t\over {2a}}}\sin^2{{2\pi J}\over t}\right| {N\over {\ell t}} 
\EssBdd \sum_{t|q,t\,\hbox{\piccolissimo LARGE}}{N\over {a\ell}}. 
$$
\par
\noindent
All these terms and the previous (with $t$ {\stampatello large}, say) give (using \thinspace $\sum_{q}{{g(\ell tq)}\over q}\EssBdd 1$)
$$
\EssBdd \sum_{T/\ell<t\le 2Q/\ell}{1\over t}\left({1\over a}+{a\over t}\right){N\over {\ell}}
\EssBdd {N\over {\ell a}} + {{Na}\over T}
$$
\par
\noindent
contributes to the {\stampatello large} $t$ in $(\ast)$. The first is diagonal-type, the others make $\EssBdd Nh^3/T$ for $S_2$.
\par
Thus, in all, the sum $S_2 \EssBdd Nh + (Q+T)T^{3/2}h + Nh^3/T$ and (here $T$ is somehow small respect to $Q$) we choose $T=(Nh^2/Q)^{2/5}$ (optimally). Also, this gives us $h=o(T)$, from $Q=o(N/\sqrt{h})$.
\par
As regards $S_1$, we may write, with $A_1(x)\defineq x$, $\forall x\in [-1/2,1/2]$ ({\stampatello odd} and) $1-${\stampatello periodicized}
$$
{\sum_{j}}'\cos{{2\pi j}\over q}{\sum_{k}}'\cos{{2\pi k\overline{q}}\over a}{\sum_{jm\equiv k(q)}}^{\ast}f(\ell m) = {\sum_{j}}'\cos{{2\pi j}\over q}\msum f(\ell m) \cos\left({{2\pi \overline{q}}\over a}qA_1\!\!\left( {{jm}\over q}\right)\right) 
$$
$$
= {\sum_{j}}'\cos{{2\pi j}\over q}\msum f(\ell m) \cos\left({{2\pi}\over a}A_1\!\!\left( {{jm}\over q}\right)\right) = {\sum_{j}}'\cos{{2\pi j}\over {aq}} \msum f(\ell m)\cos{{2\pi j\overline{m}}\over q} 
$$
\par
\noindent
where $'$ means: in $[-q/2,q/2]$ and not zero. Here $0<|k|\le q/2$ and $k\equiv jm(q)$ force $k$ to be $qA_1(jm/q)$. Then, 
$$
\sum_{0<|j|\le {q\over 2}}\cos{{2\pi j}\over {aq}} \msum f(\ell m)\cos{{2\pi j\overline{m}}\over q} \ll {1\over a}\msum |f(\ell m)|{1\over {\left \Vert {{\overline{m}}\over q}\right \Vert}}\EssBdd {1\over a}\msum {1\over {\left \Vert {{\overline{m}}\over q}\right \Vert}}
$$
\par
\noindent
using quoted [D] and this $\Rightarrow $ $S_1$ of {\stampatello diagonal-type}:
$$
S_1\EssBdd Nh.\enspace \square
$$
\bigskip

\par				
\noindent
\centerline{\bf References}
\medskip

\item{\bf [C]} \thinspace Coppola, G.\thinspace - \thinspace {\sl On the correlations, Selberg integral and symmetry of sieve functions in almost all short intervals} \thinspace - \thinspace http://arxiv.org/abs/0709.3648v3 - 9pp. (electronic).
\smallskip

\item{\bf [D]} \thinspace Davenport, H.\thinspace - \thinspace {\sl Multiplicative Number Theory} \thinspace - \thinspace Third Edition, GTM 74, Springer, New York, 2000. $\underline{{\tt MR\enspace 2001f\!:\!11001}}$
\smallskip

\item{\bf [I-K]} \thinspace Iwaniec, H. and Kowalski, E.\thinspace - \thinspace {\sl Analytic Number Theory} - American Mathematical Society Colloquium Publications, 53. AMS, Providence, RI, 2004. xii+615 pp. ISBN: 0-8218-3633-1 $\underline{{\tt MR\enspace 2005h\!:\!11005}}$
\smallskip

\item{\bf [T]} \thinspace Tenenbaum, G. \thinspace - \thinspace {\sl Introduction to analytic and probabilistic number theory} - Cambridge studies in advanced mathematics : 46, Cambridge University Press, 1995. $\underline{{\tt MR \enspace 97e\!:\!11005b}}$
\smallskip

\item{\bf [V]} \thinspace Vinogradov, I.M.\thinspace - \thinspace {\sl The Method of Trigonometrical Sums in the Theory of Numbers} - Interscience Publishers LTD, London, 1954. $\underline{{\tt MR \enspace 15,941b}}$
\medskip

\leftline{\tt Dr.Giovanni Coppola}
\leftline{\tt DIIMA - Universit\`a degli Studi di Salerno}
\leftline{\tt 84084 Fisciano (SA) - ITALY}
\leftline{\tt e-mail : gcoppola@diima.unisa.it}

\bye